\documentclass[11pt,letterpaper]{amsart}
\usepackage{amscd,times}
\usepackage[utf8]{inputenc}
\usepackage{amsmath}
\usepackage{amsthm}
\usepackage{mathtools}
\usepackage{amsfonts}
\usepackage{comment}
\usepackage{graphicx}
\usepackage{enumerate}
\usepackage{url}
\usepackage{tikz}
\usetikzlibrary{cd}
\usepackage[hidelinks]{hyperref}
\usepackage{cleveref}
\usepackage{makecell}

\usepackage[cmtip,all]{xy}
\newcommand{\longsquiggly}{\xymatrix{{}\ar@{~>}[r]&{}}}

\hypersetup {
	pdfpagemode = {UseNone},
	pdftitle = {group-related formalities},
	pdfauthor = {Giovanni Placini, Jonas Stelzig, Leopold Zoller},
	pdflang = {en-UK}
} 

\newtheorem{thmintro}{Theorem}

\newtheorem{prop}{Proposition}[subsection]
\newtheorem{lem}[prop]{Lemma}

\newtheorem{definition}[prop]{Definition}
\newtheorem{thm}[prop]{Theorem}
\newtheorem{cor}[prop]{Corollary}
\theoremstyle{remark}

\newtheorem{rem}[prop]{Remark}
\newtheorem{ex}[prop]{Example}

\newtheoremstyle{colon}
{}
{}
{\itshape}
{}
{\bfseries}
{:}
{ }
{}

\theoremstyle{colon}

\newcommand{\Q}{\mathbb{Q}}

\newcommand{\C}{\mathbb{C}}

\newcommand{\CP}{\mathbb{CP}}
\newcommand{\R}{\mathbb{R}}

\newcommand{\ft}{\mathfrak{t}}

\newcommand{\Hom}{\operatorname{Hom}}

\newcommand{\pr}{\operatorname{pr}}
\newcommand{\im}{\operatorname{im}}

\newcommand{\del}{\partial}

\newcommand{\delbar}{{\bar\partial}}

\usepackage[all]{xy}

\newcommand{\id}{\mathrm{Id}}

\newcommand{\cbba}{\operatorname{cbba}}
\newcommand{\cdga}{\operatorname{cdga}}

\newcommand{\Ho}{\operatorname{Ho}}
\newcommand{\tot}{\operatorname{to}}

\makeatother

\title{Strong formality of toric and homogeneous compact K\"ahler manifolds}
\author{G.~Placini, J.~Stelzig, L.~Zoller}

\begin{document}

\maketitle
\begin{abstract}
    All compact Kähler, or even $\partial\bar\partial$-manifolds, are rationally formal. Not all of them are strongly formal. Yet some of them are: For complete smooth complex toric varieties and homogeneous compact Kähler manifolds we show the stronger property that they are both rationally and strongly formal in a compatible way.
\end{abstract}
\section{Introduction}

All complex manifolds satisfying the $\partial\bar\partial$-Lemma, such as compact K\"ahler manifolds, are  formal in the sense of rational homotopy theory \cite{DGMS}, \cite{Sullivan}. I.e., their commutative differential graded algebra (cdga) of forms is connected by a chain of multiplicative quasi-isomorphisms to its cohomology. In particular, all Massey products vanish on such manifolds and, in the simply connected case, $\pi_{\geq 2}\otimes\Q$ can be computed from the rational cohomology ring. 

On the other hand, pluripotential homotopy theory \cite{MilSt_bigrform}, \cite{StePHT}, gives rise to a stronger notion of formality on complex manifolds, which takes into account the complex structure. 
Namely, a complex manifold is called strongly formal if its commutative, bigraded, bidifferential algebra (cbba) of smooth $\C$-valued forms is connected via real maps inducing isomorphisms in Bott-Chern and Aeppli cohomologies to a cbba with trivial differentials. 
Strong formality implies the $\del\delbar$-Lemma and (hence) formality, but not vice versa \cite{SfTo_DBC}, \cite{MilSt_bigrform}. 
On simply connected strongly formal manifolds, the real Hodge structure on $\pi_{\geq 2}\otimes \R$ can be computed from the real cohomology ring. 
Obstructions to strong formality are given by ABC-Massey products, a secondary variant of Massey products sensitive to the complex structure introduced in \cite{AnToBCform}.

It was shown in \cite{PSZ24} that compact K\"ahler (and even projective) manifolds may carry non-trivial ABC-Massey products, and hence are generally not strongly formal. On the other hand, to date only a few relatively special classes of strongly formal manifolds are known: Hermitian symmetric spaces \cite{MilSt_bigrform}, K\"ahler solvmanifolds \cite{SfToBCK} and compact K\"ahler manifolds of dimension $\geq 2$ with the Hodge diamond of a complete intersection. Since the constructions of \cite{PSZ24} are rather general, one may raise the question how restrictive a notion strong formality is and whether it is satisfied in other geometrically interesting cases. Our main results address this question as follows:

\begin{thmintro}\label{thm toric}
Complete smooth complex toric varieties are strongly formal over $\Q$.
\end{thmintro}

\begin{thmintro}\label{thm homogeneous}
Compact K\"ahler homogeneous spaces are strongly formal over $\Q$.
\end{thmintro}
We note that we do not require toric varieties to be projective, so \Cref{thm toric} includes examples of strongly formal manifolds which are not K\"ahler.

In both \Cref{thm toric} and \Cref{thm homogeneous}, the classes of spaces were previously known to be rationally formal for two reasons. Firstly, they are K\"ahler, or at least satisfy the $\partial\bar\partial$-Lemma (c.f. \cite[Cor. 5.23]{DGMS}) and are hence rationally formal. As discussed above, the $\partial\bar\partial$-Lemma is not sufficient for strong formality, so this line of reasoning breaks down in the pluripotential setting. However, in each case there also exists an independent argument for the rational formality that makes use of the group actions and their effects on the topology. We show that the analogy between pluripotential and rational homotopy theory is strong enough to adapt these equivariant arguments to the holomorphic world. In a main step of the proof of \Cref{thm toric} we argue somewhat differently from the existing equivariant proof for rational formality \cite{PanovRay}, and this line of argument would also yield an alternative proof of the rational case (cf. \Cref{rem: comparison old proof toric}).

The `over $\Q$' in the above theorems refers to a strengthening of the notion of strong formality which we introduce in this paper. Already without this addition, \Cref{thm toric} and \Cref{thm homogeneous} are new. In fact, we first prove strong formality and then show how the argument can be extended to prove strong formality over $\Q$.

Roughly speaking, strong formality over $\Q$ is formality in the category of triples $(A_\Q, A_\C, \varphi)$, where $A_\Q$ is rational cdga, $A_\C$ is a cbba and $\varphi$ is a chain of quasi-isomorphisms between $A_\Q\otimes\C$ and $A_\C$.

To every complex manifold satisfying the $\partial\bar\partial$-Lemma, we can associate two such triples: One is given by Sullivan's polynomial forms $A_{PL}(X)$, the cbba of $\C$-valued forms $A(X)$ and the de-Rham comparison quasi isomorphism. The other is given by singular cohomology, de-Rham cohomology with its Hodge-theoretic bigrading, and the comparison between the two. Strong formality over $\Q$ then means having an isomorphism in an appropriate homotopy category of such triples. Essentially by definition, there are implications
\[
\text{strong formality}/_\Q\Longrightarrow\text{strong formality}\Longrightarrow \text{rational formality}~.
\]
The converses are not true: For the second implication, this follows from the results of \cite{SfTo_DBC}, \cite{MilSt_bigrform}, or \cite{PSZ24}, since formal complex manifolds, even such which satisfy the $\partial\bar\partial$-Lemma, can have non-trivial ABC-Massey products. For the first implication, we give a counterexample in the present paper. 
\begin{thmintro}\label{thm: sfQ}
    Given a triple $(A_\Q,A_\C,\varphi)$, strong formality of $A_\C$ does not imply strong formality over $\Q$.
\end{thmintro}

Just like the manifolds in the focus of \Cref{thm toric} and \Cref{thm homogeneous}, the example we produce has only diagonal degree cohomology. The non-vanishing obstructions are computed in terms of the interplay of Aeppli and Bott-Chern cohomology with the rational structure. We show they can also be interpreted as giving rise to non-trivial extensions of rational Mixed Hodge structures on $\pi_4$ of our example.

\section{Formality and higher operations on complex manifolds}\label{sec: formality and higher ops}
\subsection{Cdga's and formality}

For a field $k$, we denote by $\cdga_k$ the category of cdga's defined over $k$ and by $\Ho(\cdga_k)=\cdga_k[quiso^{-1}]$ its homotopy category. We refer to the papers \cite{DGMS}, \cite{Sullivan}, or the textbooks \cite{FHT_RHT}, \cite{GM_RHT} for all details and background information needed.

\subsection{Bicomplexes, cbba's and strong formality}

We refer to \cite{StStrDbl} and \cite{StePHT} for details.
    A bicomplex, or double complex, $A=(A^{\bullet,\bullet},\del,\delbar)$ consists of a bigraded $\C$-vector space $A^{\bullet,\bullet}$, together with two endomorphisms $\del$ and $\delbar$ of degree $(1,0)$, resp. $(0,1)$, such that $\partial^2=\delbar^2=\del\delbar+\delbar\del=0$ or, equivalently $d^2=0$ for $d:=\del+\delbar$. 
    It is said to have a real structure if it is equipped with an anti-linear involution $\sigma:A\to A$ such that $\sigma A^{p,q}= A^{q,p}$ and $\sigma d \sigma= d$. 
    A map of bicomplexes is a linear map of the underlying vector spaces $f:A\to B$ s.t. $f(A^{p,q})\subseteq B^{p,q}$ and $fd=df$. 
    If $A,B$ are equipped with real structures $\sigma_A,\sigma_B$, $f$ is said to be real, if $f\sigma_A=\sigma_B f$. 
    A commutative, bigraded bidifferential algebra (with real structure), or for short $(\R-)$cbba, is a commutative algebra object in the category of bicomplexes (with real structure). 
    Concretely, it consists of a bicomplex $A$ (with real structure) with a (real) map $A\times A\to A$ of bicomplexes and a distinguished element $1\in A^{0,0}$ (resp. $1\in (A^{0,0})^{\sigma}$) satisfying relations of associativity, unitality and graded commutativity. 
    A (real) map of cbba's is a (real) map of bicomplexes compatible with multiplication and unit. 
    The main examples of interest for us are the $\R$-cbba's $A(X)$ of smooth, $\C$-valued differential forms on a complex manifold $X$, equipped with the conjugation action.
    Associated with any bicomplex $A$, there are various cohomology vector spaces. Of main interest for us are the graded vector space of total, or de Rham, cohomology 
    \[
    H_{dR}(A)=\frac{\ker d}{\im d}~,
    \] 
    the bigraded vector spaces of Dolbeault, conjugate Dolbeault (i.e. column and row cohomology) 
    \[
    H_{\delbar}(A):=\frac{\ker\delbar}{\im\delbar}~,\qquad H_{\del}(A)=\frac{\ker\del}{\im\delbar}~,
    \]
    and the bigraded vector spaces of Bott-Chern and Aeppli cohomology
    \[
    H_{BC}(A)=\frac{\ker\del\cap\ker\delbar}{\im\del\delbar}\qquad\mbox{and}\qquad H_A(A)=\frac{\ker\del\delbar}{\im\del+\im\delbar}~.
    \]
    These constructions are functors from bicomplexes to (bi)-graded vector spaces and we write $H_{BC}(f)$ etc.\ for the induced maps. A map of bicomplexes $f:A\to B$ is called a bigraded, or pluripotential, quasi-isomorphism if the induced maps $H_{BC}(f)$ and $H_A(f)$ are isomorphisms. The terminology pluripotential quasi-isomorphism is motivated by the potential ambiguity of the term bigraded quasi-isomorphism and the fact that Bott-Chern and Aeppli cohomology control existence and essential uniqueness of solutions to the equation $\del\delbar y=x$. If $f$ is a pluripotential quasi-isomorphism, then $H_{\del}(f)$, $H_{\delbar}(f)$ and $H_{dR}(f)$ are isomorphisms. Conversely, assuming both $A$ and $B$ have bounded antidiagonals, e.g. if they are first quadrant, then if $H_{\del}(f)$ and $H_\delbar(f)$ are isomorphisms, $f$ is a pluripotential quasi-isomorphism. In particular, if $f$ is real and $A,B$ have bounded antidiagonals, $f$ is a pluripotential quasi-isomorphism if and only if $H_\delbar(f)$ is an isomorphism. Finally, recall that all the cohomologies above are related by maps induced by the identity
    \[
    \begin{tikzcd}
        &H_{BC}(A)\ar[rd]\ar[d]\ar[ld]&\\
        H_{\delbar}(A)\ar[rd]&H_{dR}(A)\ar[d]&H_{\del}(A)\ar[ld]\\
        &H_A(A)&
    \end{tikzcd}
    \]
    If the top central map is injective, then all maps are isomorphisms \cite{DGMS}. In this case, $A$ is said to satisfy the $\partial\bar\partial$-property. If $A,B$ satisfy the $\del\delbar$-property, then a map of bicomplexes $f:A\to B$ is a pluripotential quasi-isomorphism if and only if $H_{dR}(f)$ is an isomorphism. We stress that we require $f$ to be compatible with the bigrading (so it is not true that any quasi-isomorphism between the total complexes of $A$ and $B$ is a pluirpotential quasi-isomorphism).

    A map of cbba's is called a pluripotential, or bigraded, quasi-isomorphism if its underlying map of bicomplexes is one and we can consider the homotopy category of $\R$-cbba's arising by localizing at these quasi-isomorphisms \cite{StePHT}. A cbba (resp. $\R$-cbba) $A$ is called strongly formal (over $\R$) if it isomorphic to an algebra without differentials in this homotopy category, i.e. concretely if it can be connected by a chain of (real) pluripotential quasi-isomorphisms to a (real) cbba $H$ on which $\del_H=\delbar_H=0$ \cite{MilSt_bigrform}. Without loss of generality, we may take $H=H_{BC}(A)$ with trivial differentials in this definition. A compact complex manifold $X$ is said to be strongly formal if $A(X)$ is strongly formal.
   Like for cdga's, there exists a minimal cofibrant model for any cohomologically simply connected ($\R$-)cbba $A$, which is unique up to isomorphism. I.e. there is a pluripotential quasi-isomorphism $\Lambda W\to A$ where $\Lambda W$ is a free bigraded algebra where one can choose an ordering on a basis of the generators $W$ such that $\del$ and $\delbar$ of every generator are sums of products of previous generators and such that the image of $\del\delbar$ is contained in the space of decomposable elements $\Lambda^{\geq 2}W$. The derived functor of indecomposables is called the homotopy bicomplex and can be computed as $\pi^{\ast,\ast}(A)\cong\Lambda^{\geq 1}W/\Lambda^{\geq 2}W\cong W$ or, in other words, $W$ equipped with the linear part of the differential(s).

\subsection{Strong formality over $\Q$}

\begin{definition}
    A cbba with rational structure is a triple $A=(A_\Q, A_\C, \varphi)$, where
    \begin{enumerate}
        \item $A_\Q$ is a rational cdga.
        \item $A_\C$ is an $\R$-cbba. 
        \item $\varphi$ is an isomorphism to the total cdga
        \[A_\Q\otimes\R\cong \tot(A_\C)^{\sigma=\id}\]
        in $\Ho(\cdga_\R)$.
    \end{enumerate}
    A morphism $f:(A_\Q, A_\C, \varphi)\to (B_\Q, B_\C,\psi)$ is a pair $(f_\Q, f_\C)$ where 
    \begin{enumerate}
        \item $f_\Q: A_\Q\to B_\Q$ is a morphism in $\cdga_\Q$
        \item $f_\C: A_\C\to B_\C$ is a morphism of $\R$-cbba's
    \end{enumerate}
    s.t. the diagram
    \[
    \begin{tikzcd}
        A_\Q\otimes\R\ar[d,"\varphi"]\ar[r,"f_\Q"]&B_\Q\otimes\R \ar[d,"\psi"]\\
        \tot(A_\C)^{\sigma=\id}\ar[r,"f_\C"]&  \tot(B_\C)^{\sigma=\id}
    \end{tikzcd}
    \]
    commutes in $\Ho(\cdga_\R)$.  Such a morphism is called a quasi-isomorphism if $f_\C$ is a pluripotential quasi-isomorphism (in which case also $f_\Q$ is a rational quasi-isomorphism).
    \end{definition}
\begin{rem}
    A morphism $\varphi:A\to B$ in the homotopy category of (say, connected) cdga's can be represented by a roof $A\overset{\tilde\varphi}{\longleftarrow}\Lambda V\longrightarrow B$, where the right hand map is a cofibrant (e.g. minimal) model for $B$. By abuse of notation, we will often identify $\varphi$ and $\tilde\varphi$.
\end{rem}
\begin{ex}
    To every compact complex manifold $X$ satisfying the $\partial\bar\partial$-Lemma, one associates two natural cbba's with rational structure: $(A_{PL}(X), A(X), \varphi)$ and $(H_{sing}(X,\Q), H_{dR}(X;\C), \psi)$ where $H_{dR}(X;\C)$ is considered as a bigraded algebra via the isomorphism $H_{BC}\to H_{dR}$ and equipped with trivial differential and $\varphi,\psi$ are the de Rham comparison isomorphisms.
\end{ex}

\begin{definition}
    
    A cbba with rational structure $(A_\Q, A_\C, \varphi)$ s.t. $A_\C$ satisfies the $\partial\bar\partial$-Lemma is called strongly formal over $\Q$ if it is connected by a chain of quasi-isomorphisms to its cohomology $(H(A_\Q), H_{BC}(A_\C), H(\varphi))$.
\end{definition}
\begin{definition}
    A compact complex manifold satisfying the $\partial\bar\partial$-Lemma is called strongly formal over $\Q$ if $(A_{PL}(X), A(X), \varphi)$ is strongly formal over $\Q$.
\end{definition}

\begin{definition}
    A cbba with rational structure $A=(A_\Q,A_\C,\varphi)$  is called minimal if $A_\Q$ and $A_\C$ is a minimal cdga, resp. minimal cbba.
\end{definition}
\begin{ex}\label{ex: Mimo for cbba/Q}
    Given a cbba with rational structure $A=(A_\Q, A_\C, \varphi)$, let $\alpha:M=\Lambda V\to A_\Q$ and $\beta:N=\Lambda W\to A_\C$ be rational, resp. pluripotential, minimal models. The latter exists e.g. if $A_\C$ is cohomologically simply connected\footnote{cf. \cite{StePHT}. In particular, this condition is satisfied if $A_\C$ satisfies the $\partial\bar\partial$-Lemma and $H^1=0$.}. We may choose a map $\psi:M\otimes\R\to \tot(N)^{\sigma=\id}$ such that the diagram
    \[
        \begin{tikzcd}
        \tot(N)^{\sigma=\id}\ar[d,"\beta"]&M\otimes \R\ar[l,swap,"\psi"]\ar[r,equal]\ar[d,equal]&M\otimes\R\ar[d,"\alpha"]\\
                \tot(A_\C)^{\sigma=\id}&M\otimes\R\ar[l,"\varphi"]\ar[r,swap,"\alpha"]&A_\Q\otimes\R
        \end{tikzcd}
    \]
    commutes up to homotopy. I.e. whenever $A_\Q, A_\C$ have minimal models, we can find a minimal model for $A$, i.e. a minimal cbba with rational structure which is quasi-isomorphic to $A$.
\end{ex}

\begin{prop}\label{prop: formality of mini cbba/Q}
    A minimal cbba with rational structure $(A_\Q, A_\C, \varphi)$ is strongly formal over $\Q$ if and only if the following conditions are satisfied:
    \begin{enumerate}
        \item $A_\C$ is strongly formal, i.e. there is a pluripotential quasi-isomorphism $p:A_\C\to H_{BC}(A_\C)\cong H_{dR}(A_\C)=:H(A_\C)$.
        \item $A_\Q$ is rationally formal, i.e. there is a rational quasi-isomorphism $q:A_\Q\to H(A_\Q)$.
        \item For any choice of $p,q$ as above, such that $H_{dR}(\varphi )\circ H_{dR}(q)= H_{dR}(p)\circ H_{dR}(\varphi)$ there exist an automorphism $a:A_\Q\to A_\Q$ of cdgas and an automorphism $b:A_\C\to A_\C$ of $\R$-cbbas s.t.\ $H_{dR}(a)=H_{dR}(b)=\id$ and the following diagram commutes up to homotopy:
        \[
        \begin{tikzcd}
            \tot(A_\C)^{\sigma=\id}\ar[d,"p"]&A_\Q\otimes\R\ar[l,swap,"b\circ \varphi\circ a"]\ar[d,"q"]\\
            H(A_\C)^{\sigma=\id}&H(A_\Q)\otimes\R\ar[l, "H(\varphi)"]
        \end{tikzcd}
        \]
    \end{enumerate}
\end{prop}
\begin{proof}
    That these conditions imply strong formality over $\Q$ is clear. Conversely, if $(A_\Q,A_\C,\varphi)$ is strongly formal over $\Q$, this means there exists an isomorphism in $\Ho(\R\mbox{-}\cbba)$: $p':A_\C\to H(A_\C)$ and an isomorphism in $\Ho(\cdga_\Q)$: $q':A_\Q\to H(A_\Q)$ s.t.\ the diagram 
       \[
        \begin{tikzcd}
            \tot(A_\C)^{\sigma=\id}\ar[d,"p'"]&A_\Q\otimes\R\ar[l,swap,"\varphi"]\ar[d, "q'"]\\
            H(A_\C)^{\sigma=\id}&H(A_\Q)\otimes\R\ar[l, "H(\varphi)"]
        \end{tikzcd}
        \]
        commutes. We may furthermore assume that $H(p)=H(p')$, $H(q)=H(q')$ by replacing $p'$ by $H(p)\circ H(p')^{-1}\circ p'$ and $q'$ by $H(q)\circ H(q')^{-1}\circ q'$ (where we use that $H(p)\circ H(\varphi) = H(\varphi)\circ H(q)$ to ensure the above diagram still commutes). Since $A_\C$, resp. $A_\Q$ are cofibrant, we may represent $p', q'$ by actual maps (instead of zigzags of maps) which we denote by the same name.
        By the lifting property, we find a pluripotential quasi-isomorphism $b:A_\C\to A_\C$ s.t. $p\circ b=p'$ and a rational quasi-isomorphism $a:A_\Q\to A_\Q$ s.t. $q'\circ a=q$. Since $A_\C$, $A_\Q$ are minimal, $a,b$ are actual isomorphisms and since $H(p)=H(p')$ and $H(q)=H(q')$, $a,b$ induce the identity in cohomology.   
\end{proof}

\section{Proof of \Cref{thm toric}}\label{sec: proof thm A}
We first prove strong formality. The refinement to strong formality over $\Q$ will then be developed afterwards. In principle both arguments follow the same idea although they require very different flavours of technicalities. The compatibility of both approaches then comes down to independent cohomological arguments.

\subsection{Pluripotential models for the Borel fibration}
In view of the goal of the paper, we restrict to the case of torus actions. Let $T=(S^1)^n\subseteq T_\C=(\C^\ast)^n$ and $X$ a compact complex manifold with a smooth action of $T$ through biholomorphisms (which automatically extends to an action of $T_\C$).  Let $ET\to BT$ (respectively $ET_\C\to BT_\C$) be a universal $T$-bundle (resp. $T_\C$-bundle). For concreteness, we may take $(S^\infty)^n\to (\CP^\infty)^n$, resp. $(\C^\infty-\{0\})^n\to (\CP^\infty)^n$. We write \[\pi:X_T=(X\times ET)/T\to BT\] for the Borel fibration of $X$. We will be interested in cbba models for the Borel fibration. That these should exist is plausible by the homeomorphism
\[
(X\times ET)/T\cong (X\times ET_\C)/T_\C~,
\]
where the right hand side manifestly wants to be some sort of complex space. We will not need to make this precise, but instead work directly with the cbba models.
 The so-called Cartan model, which is well-known to be a  real cdga model for the Borel construction, naturally carries the structure of a cbba and will be our model of choice. For a classic introduction see \cite{GuilleminSternberg}. The cbba structure was also used in \cite[§5.1]{Lillywhite}. 

Let $A(X)$ denote the cbba of smooth $\C$-valued forms on $X$ and $A(X)^T$ the subalgebra of $T$-invariant forms. As the $T$-action is by biholomorphisms, the action on forms respects the bigrading and $A(X)^T$ is in fact a sub-cbba. Let $\ft=Lie(T)$. As a bigraded algebra, the Cartan model is given by
\begin{equation}\label{eqn: Cartan Model}
C_T(X):=C_T(X;\C):=R\otimes_\C A(X)^T~,
\end{equation}
where $R:=S(\ft^\vee)$ denotes the symmetric algebra (with complex coefficients) on $\ft^\vee$, concentrated in degree $(1,1)$. Fixing a basis $\xi_1,...,\xi_n$ for $\ft$ and a dual basis $\xi^1,...,\xi^n$, we can identify $R=\C[\xi^1,...,\xi^n]$ and the differential $d_T$ on \eqref{eqn: Cartan Model} is determined by
\[
d_T(1\otimes \omega)=1\otimes d\omega - \sum_{i=1}^n \xi^i\otimes\iota_{\xi_i}\omega
\]
and the Leibniz rule. Note that we can write every $\xi_i=\xi_i^{1,0}+\xi_i^{0,1}$. Thus, $i_{\xi_i}$ is of type $(-1,0)+(0,-1)$ and $d_T$ is indeed of type $(1,0)+(0,1)$. E.g. the $(1,0)$-part is given by
\[
\del_T(1\otimes \omega)=1\otimes \del\omega - \sum_{i=1}^n \xi^i\otimes\iota_{(\xi_i)^{0,1}}\omega~.
\]

\begin{rem}
To express the differential without the choice of a basis, identify $C_T(X)$ with the space of polynomial maps $\ft\to A(X)^T$. Given such a map $\omega$, the differential is given by
\[
d_T(\omega)(\xi)=d(\omega(\xi))-\iota_\xi\omega(\xi)~.
\]
\end{rem}

Our next aim is to study the $\partial\bar\partial$-Lemma on the Cartan model.
The following is a standard result, going back to Cartan in the de Rham case, see e.g. \cite{CheEil_Lie}.
\begin{lem}\label{lem: invariants deldelbar}
    Let $X$ be a complex manifold satisfying the $\del\delbar$-Lemma and $G$ a connected compact Lie group acting on $X$ by biholomorphisms. Then the inclusion of the invariant forms $i:A(X)^G\hookrightarrow A(X)$ is a pluripotential quasi-isomorphism of cbba's.
\end{lem}
\begin{proof}
    Since the action is by biholomorphisms the induced action on forms respects the bigrading and there is a map of bicomplexes $a:A(X)\to A(X)^G$ given by averaging $\omega\mapsto a\omega:=\int_{G}g^*\omega dg$ with respect to the normalized Haar measure. By construction, $a\circ i=\id$, so the bicomplex $A(X)^G$ is a direct summand in $A(X)$ and hence satisfies the $\del\delbar$-Lemma as well. Now, it is known  that $a: A(X)\to A(X)^G$ induces an isomorphism in de Rham cohomology \cite{CheEil_Lie}. Since  source and target satisfy the $\partial\bar\partial$-Lemma, the averaging map is then also a pluripotential quasi-isomorphism.
\end{proof}

\begin{lem}\label{lem: dkilling}
    Let $(A,\del,\delbar)$ be a double complex satisfying the $\del\delbar$-Lemma and $a\in A^{p+1,q}\oplus A^{p,q+1}$ be $d$-exact. Then there is $\beta\in A^{p,q}$ with $d\beta = \alpha$.
\end{lem}

\begin{proof}
    Write $a=\alpha+\alpha'$ with $\alpha\in A^{p+1,q}$, $\alpha'\in A^{p,q+1}$. Then $\delbar\alpha=-\del \alpha'$ is $\del$-closed and $\delbar$-closed. By the $\del\delbar$-Lemma we find $x\in A^{p,q}$ with $\del\delbar x=\del\alpha'$. In particular $a-dx$ is $\del$-closed and $\delbar$-closed. As it is again $d$-exact it is also $\del\delbar$-exact and we find $y_1,y_2$ with $\del\delbar y_1=\alpha-\del x$ and $\del\delbar y_2= \alpha'-\delbar x$. But then $\beta = \delbar y_1 - \del y_2+x$ does the job.
\end{proof}

\begin{lem}\label{lem: BCsurjectivity}
    Let $T$ act by biholomorphisms on a compact complex manifold $X$ satisfying the $\partial\bar\partial$-Lemma. If $C_T(X)\rightarrow A(X)^T$ induces a surjection in de Rham cohomology, then it induces a surjection in Bott-Chern cohomology.
\end{lem}

\begin{proof}
    Let $\theta\in (A(X)^T)^{p,q}$ be a $\del$-closed and $\bar\del$-closed form inducing a non-trivial class in $H(X)$. We argue that it can be extended to a $\del$- and $\delbar$-closed equivariant form in $C_T(X)^{p,q}$.

    Assume inductively that for some $k\geq 1$ we have found $\theta_i\in R^{2i}\otimes (A(X)^T)^{p-i,q-i}$, $1\leq i\leq k-1$ such that for $\tilde \theta = \theta +\theta_1+\ldots+\theta_{k-1}$ we have $d_T(\tilde \theta)\in R^{2k}\otimes A(X)^T$. We write
    \[d_T(\tilde\theta)=\sum f_i (\eta_i+\eta_i')~,\]
    where the $f_i$ are a basis of $R^{k,k}$, $\eta_i\in (A(X)^T)^{p-k+1,q-k}$, and $\eta_i'\in (A(X)^T)^{p-k,q-k+1)}$.
    Cohomological surjectivity of $C_T(X)\rightarrow A(X)^T$ is equivalent to the collapse of the Serre spectral sequence of the Borel fibration at the $E_2$-page, i.e.\ of the spectral sequence that arises from the decreasing filtration $R^{\geq *}\otimes A(X)^T$. In particular it follows that $d_T(\tilde\theta)$ and hence all the individual $\eta_i+\eta_i'$ are $d$-exact as otherwise the spectral sequence would have a non-trivial differential on the $E_{2k}$-page. By Lemma \Cref{lem: invariants deldelbar} $A(X)^T$ satisfies the $\del\delbar$-Lemma and hence by \Cref{lem: dkilling} we find $\beta_i\in (A(X)^T)^{p-k,q-k}$ with $d\beta_i=\eta_i+\eta_i'$. For $\theta_k=\sum -f_i\beta_i$ we have $d_T(\tilde\theta+\theta_k)\in R^{2(k+1)}\otimes (A(X)^T)^{p,q}$.

    By the induction we obtain an extension of $\theta$ to a $d$-closed equivariant form in $C_T(X)^{p,q}$ which is then automatically $\del$-closed and $\delbar$-closed.\end{proof}

\begin{prop}\label{prop: catan deldelbar}
    Let $T$ act by biholomorphisms on a compact complex manifold $X$ satisfying the $\partial\bar\partial$-Lemma. If $C_T(X)\rightarrow A(X)^T$ induces a surjection in de Rham cohomology, then $C_T(X)$ satisfies the $\del\delbar$-Lemma.
\end{prop}

\begin{proof}
    By \Cref{lem: BCsurjectivity} we find $\del$-closed and $\delbar$-closed forms $\theta_1,\ldots,\theta_k\in C_T(X)$ whose restrictions to $A(X)^T$ induce a basis of $H_{BC}(X)$. We consider the double complex $K=R\otimes \langle\theta_1,\ldots,\theta_k\rangle$ with the trivial differential and the inclusion $\iota\colon K\rightarrow C_T(X)$. The maps $H_\delbar(X)\leftarrow H_{BC}(X)\rightarrow H_\del(X)$ are isomorphisms by the $\del\delbar$-Lemma. Hence the $\theta_i$ restrict to bases for both $H_\delbar(X)$ and $H_\del(X)$. By a Leray-Hirsch-type spectral sequence argument it follows that $H_\delbar(C_T(X))$ and $H_\del(C_T(X))$ are free $R$-modules with basis given by the classes of the $\theta_i$. In particular $\iota$ induces an isomorphism in both Dolbeaut and anti Dolbeaut cohomology. Since $K$ and $C_T(X)$ are both first quadrant double complexes and $K$ satisfies the $\del\delbar$-Lemma, it follows that $C_T(X)$ satisfies the $\del\delbar$-Lemma as well.\footnote{This follows using the characterization of the $\partial\bar\partial$-Lemma as `degenerate Fr\"olicher spectral sequence(s) + opposed filtrations on de Rham' \cite{DGMS} since the map $K\to C_T(X)$ must induce a bifiltered isomorphism in de Rham cohomology as it induces isomorphisms on the associated graded for row and column filtration, or using \cite[Prop. 11]{StStrDbl} or \cite[Thm C]{StePHT} and the characterization of the $\partial\bar\partial$-Lemma as `only dots and squares' \cite{DGMS}.}
\end{proof}

\begin{rem}
    It follows from \cite{Bla_VarAnCplx} that any action on a simply connected space whose cohomology satisfies the Hard Lefschetz property also satisfies the surjectivity condition in \Cref{prop: catan deldelbar}. In particular \Cref{prop: catan deldelbar} shows that the Cartan model of a holomorphic torus action on a compact simply connected Kähler manifold satisfies the $\del\delbar$-Lemma. This also follows from a claim in the proof of \cite[Thm. 5.1.]{Lillywhite}, where only an argument for $E_1$-degeneration is given however.
\end{rem}

\subsection{Strong formality}
Now let $X$ be a complete smooth complex toric variety.
We recall some known results on the cohomology rings of toric varieties. All coefficients are complex.

The equivariant cohomology of $X$ is $H_T(X):=H(X_T)$. By construction, it is an algebra over $H(BT)=\C[t_1,...,t_n]=:R$. Denote by $D_1,...,D_m$ the torus-invariant divisors on $X$. They have fundamental classes $\tau_i\in H_T^2(X)$, defined as the image of $1\in H^0_T(D_i)$ under the pushforward $H^0_T(D_i)\to H^2_T(X)$. We call these the Thom classes of the $D_i$.

One has (\cite{DJ}, \cite{MasudaPanov}, see also \cite[Thm. 7.4.33]{BuchPan}) the following:
\begin{prop} \label{prop: toric} Let $(X,T)$ and $D_i, \tau_i$ be as above.
\begin{enumerate}
    \item Let $\mathcal{R}=\{I\subseteq\{1,...,m\}\mid \bigcap_{i\in I} D_i=\emptyset\}$. Then
    \[H_T^\ast(X)=\C[\tau_1,...,\tau_m]\bigg/\left\langle\prod_{i\in I} \tau_i\mid I\in \mathcal{R}\right\rangle~.\]
    \item The Serre spectral sequence for $\pi:X_T\to BT$ degenerates. In particular, $H_T^\ast(X)$ is a free $H^\ast(BT)$-module
  
    and restriction induces a multiplicative isomorphism
    \[H^\ast(X)=H_T^\ast(X)/(R^{>0}\cdot H_T^\ast(X))~.\]
\end{enumerate}
\end{prop}

\begin{lem}
One may choose real forms $\alpha_i\in C_T(X)^{1,1}$ such that $\tau_i=[\alpha_i]$ which are concentrated arbitrarily closely along $D_i$. In particular, one may assume $\prod_{i\in I} \alpha_i=0$ for $I\in \mathcal{R}$.
\end{lem}
\begin{proof}

The class $\tau_i$ is the first Chern class of the line bundle $L\rightarrow X$ associated to a $T$-invariant divisor $D$. We note first that the $T$-action lifts to an action on $L$. The latter is defined by a collection of local holomorphic functions $f_U\colon U\rightarrow \C$ such that $D\cap U$ is the set of zeroes of $f_U$, all of which are of first order. Note that for $t\in T$, the function $t^*f_U$ on $t^{-1}U$ shares the same properties. Hence we may extend the family of defining functions to be closed under the $T$-action. Now $L$ is defined as the line bundle with transition functions $\frac{f_U}{f_V}$ and the transformations $t\times \id_\C\colon U\times \C\rightarrow tU\times \C$ are compatible with the transition maps, hence lifting the $T$-action. We take one of the neighbourhoods for the construction to be $X\backslash D$, setting $f_{X\backslash D}=1$. In particular, on the resulting $T$-invariant bundle chart $(X\backslash D)\times \C$ the group operates only on the first component.

To construct the equivariant Chern form, we first choose a $T$-invariant hermitian metric which is the standard metric in the bundle chart $(X\backslash D)\times\C$ outside a small tube around $D$. Now we fix the Chern connection on $L$ with respect to that metric. In \cite{BerlineVergne} representatives for equivariant characteristic classes are constructed via equivariant Chern-Weyl theory. The construction considers the connection $1$-form $\omega$ on $\mathrm{Fr}(L)$, the principal $U(1)=S^1$-bundle of norm-$1$-vectors associated with $L$, which in our case takes values in $\mathfrak{u}(1)\cong \R\subseteq \C\cong\mathfrak{gl}(1,\C)$ and is $T$-invariant. Now consider the equivariant $2$-form
\[J-\frac{\Omega}{2\pi i},\]
where $\Omega$ is the usual non-equivariant curvature $2$-form of $\omega$ and $J\in \mathfrak{t^*}\otimes A^0(X)^T\otimes \mathfrak{gl}(1,\C)$ is given by sending $x\in \mathfrak{t}$ to $\omega(\overline x)$. Here $\overline{x}$ denotes the fundamental vector field on the frame bundle $\mathrm{Fr}(L)$. Equivariant characteristic forms are then obtained by inserting the equivariant curvature into invariant polynomials on the Lie algebra of the structure group -- in the case of the first Chern class this yields the element $J-\frac{\Omega}{2\pi i}\in C_T^{1,1}(\mathrm{Fr}(L))$ under the canonical identification $\mathfrak{gl}(1,\C)\cong \C$. This form is basic and induces a closed real element in $C^{1,1}_T(X)$ by \cite[Lemma 2.12, Prop. 2.13]{BerlineVergne}. To verify that this is of type $(1,1)$ we note that this concerns only the part of polynomial $\mathfrak{t}^*$-degree $0$, i.e.\ $-\Omega/{2\pi i}$ where the construction agrees with the standard construction for non-equivariant Chern forms.

Since the metric is flat when sufficiently far away from $D$, it follows that $\Omega$ is supported near $D$ and it remains to prove the same for $J$. When expressed in the bundle chart $(X\backslash D)\times \C$, the orbits of the $T$-action are constant in the second component. Since by construction the connection is just given by the exterior derivative $d$ when sufficiently far away from $D$, we deduce that those orbits are horizontal. In particular, the connection form $\omega$ vanishes on the fundamental vector fields away from $D$ and hence $J$ is supported near $D$.\end{proof}

\begin{thm}\label{thm main thm complex}
    A complete smooth toric variety is strongly formal.
\end{thm}

\begin{proof}
    There is a direct map of $\R$-cbba's from $f:H_T^{\ast,\ast}(X;\C)\to C_T(X)$, namely sending $\tau_i$ to $\alpha_i$, chosen as in the previous Lemma, ensuring well-definedness. This is a pluripotential quasi-isomorphism by the formula for $H_T(X)$ in \Cref{prop: toric}.

    We consider the composition $g\colon R\rightarrow H_T^{\ast,\ast}(X;\C)\to C_T(X)$ (which is not equal to the standard inclusion). Let $S=\Lambda (s_i,\del s_i,\delbar s_i, i=1,\ldots,n)$. Now $f$ induces a quasi-isomorphism between the extensions
    $$H_T^{\ast,\ast}(X;\C)\otimes S\to C_T(X)\otimes S~,$$
    where $i\del\delbar s_i$ is the image of $t_i$ in the respective cbba (using $g(t_i)$ on the right hand side). The fact that $H_T^{\ast,\ast}(X;\C)$ is free over $R$ implies that the projection $$H_T^{\ast,\ast}(X;\C)\otimes S\to H_T^{\ast,\ast}(X;\C)/(R^+\cdot H_T^{\ast,\ast}(X;\C))$$
    is a pluripotential quasi-isomorphism. To see this observe e.g. that $H_T^{\ast,\ast}(X;\C)\otimes S$ splits additively into a direct sum of sub-bicomplexes that are up to degree shift of the form $R\otimes S$, $i\del\delbar s_i=t_i$, and hence quasi-isomorphic to a copy of $\C$. In particular $C_T(X)\otimes S$ is formal.

    Note that by definition $g$ induces the same map as the standard inclusion $R\rightarrow C_T(X)$ on the level of cohomology. Hence $g(t_i)-t_i= \del_T\delbar_T\eta_i$ for some $\eta_i\in A(X)^{0,0}$. Now consider the automorphism of $\varphi\colon C_T(X)\otimes S$ defined by being the identity on $C_T(X)$ and mapping $s_i\mapsto s_i-\eta_i$. The original $C_T(X)\otimes S$ is isomorphic to $(C_T(X)\otimes S,\del',\delbar')$ with $\del'=\varphi^{-1}\del\varphi$, $\delbar'=\varphi^{-1}\delbar\varphi$. The latter now has $i\del'\delbar'(s_i)=t_i$ thus has the form of an extension $$C_T(X)\otimes S= (R\otimes S)\otimes A(X)$$
    of the contractible cbba $R\otimes S$, hence the projection to $A(X)$ is a quasi-isomorphism. This shows formality of $A(X)$.
\end{proof}

\begin{rem}\label{rem: comparison old proof toric}
    In \cite{PanovRay} the authors give a proof of (rational) formality of toric and quasi-toric manifolds as follows: In a first step, they  show the formality of the Borel construction using Davis-Januszkiewicz complexes. In a second step, they deduce formality of the manifold itself using degeneracy of the Leray spectral sequence associated with the Borel construction. Our proof runs along the same lines, but we give a different argument for (strong) formality of the Borel construction by building a direct map from the cohomology. With some simplifications, this method can also be applied to obtain an alternative proof of the first step in the proof of \cite{PanovRay}.
\end{rem}

\subsection{Strong formality over $\Q$}

In the proof of \Cref{thm main thm complex} we built a map \[\varphi\colon H^{*,*}(X;\C)\rightarrow A(X)^{*,*}\] to show bigraded formality of a toric manifold $X$. Domain and target of the map come with canonical rational structures in forms of the maps \[H(X;\Q)\rightarrow H(X;\C)\quad A_{PL}(X;\Q)\rightarrow A_{PL}(X;\C)\simeq A(X)\]
where the right hand quasi-isomorphism is fixed and elaborated on below (here and from now on, we will omit writing the totalization explicitly). We will show
\begin{thm}\label{thm: ratioplex toric}
The map $\varphi$ and the canonical rational structures give rise to a commutative diagram in $\Ho(\cdga_\R)$
\[
\begin{tikzcd}
    A(X)^{\sigma=\id} & \ar[l] A_{PL}(X;\Q)\otimes\R\\
    \ar[u, "\varphi"] H(X;\R) & \ar[l,"="]\ar[u] H(X;\Q)\otimes\R
\end{tikzcd}
\]

\end{thm}

To prove this we will essentially mimic the construction of $\varphi$ on rational polynomial forms and then go through comparison isomorphisms between singular and de Rham cohomology in order to certify that, up to homotopy, this gives the original $\varphi$ after tensoring with $\C$.

We need to recall some machinery following \cite[chapter 10]{FHT_RHT}. Let $k$ be some ground field of characteristic $0$. Then there are simplicial cochain algebras $A_{PL}$ and $C_{PL}$ and given any simplicial set $K$ we may form the dgas
$A_{PL}(K)$ and $C_{PL}(K)$ as the set of all simplicial morphisms from $K$ into the respective target algebra. Now $A_{PL}(K)$ is commutative and $C_{PL}(K)$ is the singular cochain algebra with coefficients in $k$. As shown in \cite[Thm. 10.9]{FHT_RHT} for any simplicial set $K$, there are isomorphisms in $\Ho(\operatorname{dga}_k)$ \[A_{PL}(K)\cong C_{PL}(K)~,\]natural in $K$.

The key to the construction of $\varphi$ was to consider representatives of Thom classes of divisors with small supports. To mirror this on the simplicial side we need to restrict to small simplices. Hence we consider simplices which respect to a fine open cover of $X$. Cover the orbit space $X/T$ (in the projective case we may think of the moment polytope) with an open cover $\overline{\mathcal{U}}$ satisfying the following property: For any facets  $F,G$ with $F\cap G=\emptyset$ and any $U\in\overline{\mathcal{U}}$ one has $F\cap U=\emptyset$ or $G\cap U=\emptyset$. Considering the preimages in $X$ gives an open cover $\mathcal{U}$ of $X$ consisting of $T$-invariant sets.

Let $S(X)$ denote the singular simplices on $X$ and $S_\mathcal{U}(X)$ the simplices subordinate to the open cover $\mathcal{U}$. Furthermore we set $A_{PL,\mathcal{U}}(X)=A_{PL}(S_\mathcal{U}(X))$. Then the inclusion $S_\mathcal{U}(X)\rightarrow S(X)$ induces a cdga morphism
\[A_{PL}(X)\rightarrow A_{PL,\mathcal{U}}(X)\]
which is a quasi-isomorphism because this holds for the induced map \[C_{PL}(S(X))\rightarrow C_{PL}(S_\mathcal{U}(X))\] on singular cohomology by the classic barycentric subdivision argument. As the cover is by $T$-invariant sets we obtain a cover of $X_T$ and a quasi-isomorphism
$A_{PL}(X_T)\rightarrow A_{PL,\mathcal{U}}(X_T)$ using the analogous notation and reasoning. Recall that $\varphi$ was constructed as a composition
\[H(X;\C)\xrightarrow{\varphi_T} C_T(X)\rightarrow A(X)~.\]
We begin by constructing the rational analogue of $\varphi_T$.

\begin{prop}
    There is a cdga morphism $\psi\colon H(X_T;\Q)\rightarrow A_{PL,\mathcal{U}}(X_T;\Q)$ inducing the identity on cohomology.
\end{prop}

\begin{proof}
    As in the proof of \Cref{thm main thm complex} we show that the Thom class of a divisor corresponding to a facet has a representative such that the supports belonging non-intersecting facets do not intersect. Then the claim follows as in the proof of \Cref{thm main thm complex} by using that $H(X_T)$ is generated by the Thom classes with the only relations coming from the intersection relations of the facets as described in \Cref{prop: toric}.

   Let $F$ be a divisor corresponding to a facet. We consider the cdga $A_{PL}(X, X-F)$ of all forms vanishing on simplices in $S(X-F)$ and use the analogous definition for $A_{PL,\mathcal{U}}(X,X-F)$. We have a commutative diagram
    \[
\begin{tikzcd}
    A_{PL}(X,X-F) \ar[r]\ar[d] & A_{PL}(X) \ar[d]\ar[r] & A_{PL}(X-F)\ar[d]\\
    A_{PL,\mathcal{U}}(X,X-F)\ar[r] & A_{PL,\mathcal{U}}(X) \ar[r]& A_{PL,\mathcal{U}}(X-F)
\end{tikzcd}
\]
As argued previously  the middle and right hand vertical arrows are quasi-isomorphisms so by the $5$-Lemma the same holds for the left hand morphism. Furthermore the Thom class of $F$ lies in the image of $H(X,X-F)\rightarrow H(X)$ and thus we conclude that in $A_{PL,\mathcal{U}}(X)$ it has a representative coming from $A_{PL,\mathcal{U}}(X,X-F)$. Choose such a representative $\tau_F$.

Now if $H$ is another facet not intersecting $F$, then from the construction of $\mathcal{U}$ it follows that every simplex in $S_\mathcal{U}(X)$ lies in $X-F$ or in $X-H$. It follows that $\tau_F\cdot\tau_H=0$. This concludes the proof of the proposition.
\end{proof}

To compare this to de Rham forms assume now that $k\supset \R$ and recall from [RHT, 11(c)] that the differential forms on the standard simplex yield a simplicial cochain algebra $A_{DR}$ which contains $A_{PL}$. It is furthermore shown that for every simplicial set $K$ the inclusion
\[A_{PL}(K)\rightarrow A_{DR}(K)\]
is a quasi-isomorphism \cite[proof of Thm 11.4 part(i)]{FHT_RHT}. Furthermore it is shown that for the set $S_\infty(X)$ of smooth simplices in $X$, the canonical restriction of differential forms to simplices
\[A(X)\rightarrow A_{DR}(S_\infty(X))\]
is a quasi-isomorphism. Setting $S_{\infty,\mathcal{U}}(X)=S_\mathcal{U}(X)\cap S_\infty(X)$ we obtain the diagram

\[
\begin{tikzcd}
    A(X) \ar[r]\ar[dr]& A_{DR}(S_\infty(X))\ar[d]\\
    & A_{DR}(S_{\infty,\mathcal{U}}(X))
\end{tikzcd}
\]

On cohomology, the vertical map agrees with the map induced by $S_{\infty,\mathcal{U}}\rightarrow S_\infty$ on singular cohomology by the aforementioned comparison isomorphisms between $A_{DR}, A_{PL}$, and $C_{PL}$. On singular cohomology it is an isomorphism by the barycentric subdivision argument. We conclude that the diagonal map in the above diagram is a quasi-isomorphism.

\begin{lem}\label{lem: homotopycom}
   The diagram 

\begin{tikzcd}
    H(X_T;\Q)\ar[r, "\psi"]\ar[d] & A_{PL,\mathcal{U}}(X_T;\Q)\ar[r] & A_{PL,\mathcal{U}}(X;\Q)\ar[dr]&\\
    H(X_T;\R)\ar[r, "\varphi_T"] & C_T(X;\R)\ar[r] & A(X;\R)\ar[r] & A_{DR}(S_{\infty,\mathcal{U}}(X;\R))
\end{tikzcd}
commutes up to homotopy.
\end{lem}

\begin{proof}
    Let $f,g$ denote the two compositions and $B:= A_{DR}(S_{\infty,\mathcal{U}}(X))$. Let $x$ be the Thom class of a facet $F$. Since $f,g$ cover the same map on the level of cohomology one has $f(x)=g(x)+dy$, where $f(x),g(x)$ are cocycles. Hence one may start defining the desired homotopy
    $$h\colon H(X_T;\Q)\rightarrow B\otimes \Lambda(t,dt)=:B(t,dt)$$
    via $x\mapsto f(x)-t(f(x)-g(x))-dt\cdot y$ which commutes with differentials and would yield a homotopy from $f$ to $g$ if it could be defined like this on all Thom classes of facets. Note that by construction $f(x)$ and $g(x)$ vanish on simplices which lie outside of a tubular neighbourhood $U$ of the divisor associated to $F$. Hence, if the same were true for $y$, then $h(x)$ would be supported near $F$ and we could define $h$ on all of $H(X_T;\Q)$ whose only relations between the Thom classes come from empty intersections of facets. It remains to prove that $y$ can be chosen such that it vanishes on simplices contained in $X-U$. However by construction $f(x),g(x)$ both represent the Thom class of $F$ in \[H^2(X,X-U;\R)\cong H^2(A_{DR}(S_{\infty,\mathcal{U}}(X,X-U);\R))\] so $y$ can be chosen in the subalgebra $A_{DR}(S_{\infty,\mathcal{U}}(X,X-U);\R)$.    
\end{proof}

Now we choose a basis $x_i$ for the image of $H^2(BT;\Q)\rightarrow H^2(X_T;\Q)$ and form $H(X_T;\Q)\otimes S$, where $S=\Lambda (s_i)$ is an exterior algebra with $ds_i = x_i$. Since the $x_i$ become exact in $A_{DR}(S_{\infty,\mathcal{U}}(X;\R))$ we can extend the diagram from \Cref{lem: homotopycom} to the first two rows in the diagram

\[
\begin{tikzcd}
    H(X_T;\Q)\otimes S\ar[d]\ar[drr, "f"] & &\\
    H(X_T;\R)\otimes S\ar[r]\ar[d]& A(X;\R)\ar[r]& A_{DR}(S_{\infty,\mathcal{U}}(X;\R))\\
    H(X;\R)\ar[ur, "\varphi"] & &
\end{tikzcd}
\]
where $f(s_i)$ is defined as some choice of $y_i$ with $d(y_i)=f(x_i)$. The triangle involving the middle and lower row is the one from the construction of $\varphi$ as in the proof of \Cref{thm main thm complex} and is in particular commutative.

\begin{lem}
    The upper triangle commutes up to homotopy.
\end{lem}

\begin{proof}
    We wish to extend the homotopy $h\colon H(X_T;\Q)\rightarrow A_{DR}(S_{\infty,\mathcal{U}}(X;\R))(t,dt)$ from \Cref{lem: homotopycom}. Let $g$ denote the composition through the middle row of the above diagram. Then one has $$h(x_i)= f(x_i)-t(f(x_i)-g(x_i))-dt\cdot y_i$$
    for some $y_i$ with $dy_i=f(x_i)-g(x_i)$. In particular $f(s_i)-g(s_i)-y_i$ is closed. Since $H^1(X;\C)=0$ it is of the form $dz_i$ for some $z_i$ of degree $0$. We set
    $$h(s_i)=f(s_i)-t(f(s_i)-g(s_i))-dt\cdot z_i$$
    which defines the desired extension.
\end{proof}

Now, since $H(X_T;\Q)\otimes S\simeq H(X;\Q)$ and $f$ factors through $A_{PL,\mathcal{U}}(X;\Q)$, after inverting $A(X;\R)\rightarrow A_{DR}(S_{\infty,\mathcal{U}}(X;\R))$, we obtain a commutative square
\[\begin{tikzcd}
    A(X;\R) & \ar[l] A_{PL,\mathcal{U}}(X;\Q)\otimes\R\\
    \ar[u, "\varphi"] H(X;\R) & \ar[l]\ar[u] H(X;\Q)\otimes\R
\end{tikzcd}\]
in $\Ho(\cdga_\R)$. The fact that one can replace the top right corner to end up with the desired diagram from \Cref{thm: ratioplex toric} follows from the commutative diagram

\[
\begin{tikzcd}
    A(X;\R) \ar[r]\ar[dr]& A_{DR}(S_\infty(X);\R)\ar[d] & \ar[l] A_{PL}(X;\Q)\ar[d]\\
     & A_{DR}(S_{\infty,\mathcal{U}}(X);\R) & \ar[l] A_{PL,\mathcal{U}}(X;\Q)
\end{tikzcd}
\]
of cdgas, where the left hand triangle and the right hand vertical map are quasi-isomorphisms.

\section{Proof of \Cref{thm homogeneous}}\label{proof thm B}

The proof of \Cref{thm homogeneous} will be a consequence of a more general intrinsic formality result for cbba's (over $\Q$) which have cohomology rings of complete intersection type, which can be seen as a pluripotential version of the main result of \cite{Body_regular}, when read in conjunction with \cite{Sullivan}.
\subsection{Intrinsic formality of algebras of complete intersection type}
Let $k$ be a field of characteristic zero.
\begin{definition}\,
    \begin{enumerate}
        \item A graded $k$-algebra $H^{\ast}$ is said to be of complete intersection type if we can write $H=k[x_1,...,x_n]/(r_1,...,r_m)$ with $x_i$ of even positive degree and $r_1,...,r_m\in k[x_1,...,x_n]$ a regular sequence.        
        \item A bigraded $k$-algebra $H^{\ast,\ast}$ is said to be of complete intersection type if we can write $H\cong k[x_1,...,x_n]/(r_1,...,r_m)$ with $|x_i|=(p_i,p_i)$ for some $p_i\in \mathbb{N}_{\geq 1}$ and $r_1,...,r_m\in k[x_1,...,x_n]$ a regular sequence consisting of elements of pure bidegree.
    \end{enumerate}
\end{definition}

The following standard result on regular sequences and their Koszul complexes (see e.g. \cite{Eisenbud_CA}) shows how to calculate minimal models of graded algebras of complete intersection type.

\begin{lem}\label{lem: singly graded Koszul}
Assume $H$ is a graded algebra of complete intersection type with a presentation $H=k[x_1,...,x_n]/(r_1,...,r_m)$.
\begin{enumerate}
    \item Let $V$ be a graded vector space with basis $p_1,...,p_m$ with $|p_i|=|r_i|-1$. Consider the cdga $K:=k[x_1,...,x_n]\otimes\Lambda V$ with differential $d$ defined by $dx_i=0$ and $dp_i=r_i$ and the Leibniz-rule. Consider $H$ as a cdga with trivial differential. Then $K\to H$ defined by $x_i\mapsto x_i$ and $p_i\mapsto 0$ is a quasi-isomorphism of cdga's.

    \item Let $A$ be a cdga with cohomology $H$. Then $A$ is formal.
\end{enumerate}
\end{lem}
We deduce a variant for bigraded algebras:
\begin{lem}\label{lem: bigraded Koszul}
    Let $H$ be a bigraded algebra of complete intersection type with presentation 
    \[
     H=\C[X_1,...,X_n]/(R_1,...,R_m)
    \]
    \begin{enumerate}
        \item Let $W$ be a graded vector space with basis $P_j,\del P_j,\delbar P_j$, $j=1,...,m$ s.t. $|P_j|=|R_j|-(1,1)$. Consider the cbba $K:=\C[X_1,...,X_n]\otimes\Lambda V$ with differential defined by the Leibniz-rule and
    \begin{align*}
        d X_i&=0\\
        d P_j&=\del P_j+\delbar P_j\\
        i\del\delbar P_j &=R_i
    \end{align*}
    Consider $H$ as a cbba with trivial differential.
    Then, the map $K\to H$, defined by $X_i\mapsto X_i$ and $P_j\mapsto 0$ is a pluripotential quasi-isomorphism.

    \item  Let $A$ be a cbba satisfying the $\del\delbar$-Lemma and with cohomology $H$. Then $A$ is strongly formal.
    \end{enumerate}
    
\end{lem}
\begin{rem}
    Consider $\C[X_1,...,X_n]=\R[X_1,...,X_n]\otimes\C$ as an $\R$-cbba with involution acting only on the coefficients. Thus, if the $R_j$ are real, $H$ has a natural structure of an $\R$-cbba and the map $K\to H$ is a map of $\R$-cbba's if $V$ carries a natural antilinear involution fixing $P_j$ and interchanging $\del R_i$ and $\delbar R_i$. 
\end{rem}
\begin{rem}
    Note that this implies in particular that $K$ satisfies the $\del\delbar$-Lemma.
\end{rem}
\begin{proof}[Proof of \Cref{lem: bigraded Koszul}]
    With regards to $(1)$ we note that both $H$ and $K$ are first quadrant cbba's. As such, we may test that $\varphi:K\to H$ is a pluripotential quasi-isomorphism by checking that $H_\delbar(\varphi)$ and $H_\del(\varphi)$ are isomorphisms. Let $U_\delbar$ denote the forgetful functor that forgets the $(1,0)$ differential, i.e. it maps a cbba $(A^{\ast,\ast},\del,\delbar, \wedge)$ to the bigraded cdga $(A^{\ast,\ast},\delbar,\wedge)$. Then $U_{\delbar}K$ can be decomposed as a tensor product $K'\otimes C$ where $K'=H\otimes V'$ with $V'=\langle \del p_j\rangle_{j=1,...,m}$ and $C=\Lambda V''$ with $V''=\langle p_j, \del p_j\rangle_{j=1,...,m}$. Now, $C$ is contractible and $U_\delbar(\varphi)|_{C}=0$, while $U_\delbar (\varphi)|_{K'}$ is a quasi-isomorphism (w.r.t $\delbar$) by \Cref{lem: singly graded Koszul}. The analogous argument shows that $H_{\del}(\varphi)$ is also a quasi-isomorphism.

    For the proof of $(2)$ we construct a pluripotential quasi-isomorphism $\Phi\colon K\rightarrow A$ as follows: We define $\Phi(X_i)$ to be equal to a diagonal dergree representative for the respective cohomology class. The image $\Phi(R_i)$ is $i\del\delbar$-exact in $A$. Now we define $\Phi(P)$ as any such $i\del\delbar $-primitive. It follows from part $(1)$ that this is a pluripotential quasi-isomorphism.
\end{proof}

Now, let us adapt everything to cbba's with rational structures:

\begin{definition}
    A bigraded algebra of complete intersection type with rational structure is a triple $(H_\Q,H_\C,\iota)$, where $H_\Q$ is a graded $k$-algebra, $H_\C$ is a bigraded $k$-algebra of complete intersection type with an antilinear involution $\sigma$ and $\iota:H_\Q\to H_\C$ is a map of algebras inducing an isomorphism $H_\Q\otimes\R\cong H_\C^{\sigma=\id}$.
\end{definition}
\begin{rem}\label{rem: cit}
    The property of being an algebra of complete intersection type is invariant under field extensions, hence also $H_\Q$ is of complete intersection type.
\end{rem}

\begin{thm}
Let $A=(A_\Q, A_\C,\varphi)$ be a cbba with rational structure such that $A_\C$ satisfies the $\partial\bar\partial$-Lemma such that $H^{\ast,\ast}_{dR}(A_\C)$ is a bigraded algebra of complete intersection type. Then $A$ is strongly formal over $\Q$.    
\end{thm}

\begin{proof}
Write $H_\Q=H(A_\Q)$ and $H_\C=H(A_\C)$. By \Cref{rem: cit} the rational cohomology is also of complete intersection type and we may choose presentations
 \[
    H_\Q=\Q[x_1,...,x_n]/(r_1,...,r_m)\quad\text{and} \quad H_\C=\C[X_1,...,X_n]/(R_1,...,R_m)
    \]
    with $H(\varphi)(x_i)=X_i$ and $H(\varphi)(r_i)=R_i$ (where the latter by abuse of notation refers to the induced map on $\Q[x_1,\ldots,x_n]$.
 We construct a rational minimal model of $K_\Q\to A_\Q$ and a pluripotential minimal model $K_\C\to A_\C$ as in Lemmas \ref{lem: singly graded Koszul} and \ref{lem: bigraded Koszul}. Now the latter restricts to a quasi-isomorphism of $\R$-cdgas on the involution fixed points, so by lifting we obtain a morphism $K_\Q\to K_\C$ of cdgas which makes the diagram

    \[
    \begin{tikzcd}   K_\Q\otimes \mathbb{R} \ar[d]\ar[r]&\tot(K_\C)^{\sigma=\id}\ar[d]\\
        A_\Q\otimes \mathbb{R} \ar[r]&\tot(A_\C)^{\sigma=\id}
    \end{tikzcd}
\]
commute in $\Ho(\cdga_\R)$. Now we define maps $K_\Q\rightarrow H_\Q$ and $K_\C\to H_\C$ as in Lemmas \ref{lem: singly graded Koszul} and \ref{lem: bigraded Koszul}. Note that the resulting square
\[\begin{tikzcd}
    K_\Q\ar[d]\ar[r] & K_\C\ar[d]\\ H_\Q\ar[r] & H_\C
\end{tikzcd}\]
commutes on the $x_i\in K_\Q$ by construction. On odd degree generators the compositions in the diagram are trivial since $H_\C^{odd}=0$. Hence the entire diagram commutes.   \end{proof}

\begin{rem}
    One can state the condition that $A$ satisfies the $\partial\bar\partial$-Lemma purely in terms of the cohomology of $A$ (namely, it is equivalent to the statement that $h_{BC}+h_A=2b_k$, \cite{AnTo}). Thus, this is a statement about intrinsic strong formality of cbba's with cohomology of complete intersection type. 
\end{rem}
\begin{cor}[=\Cref{thm homogeneous}]
    Any compact homogeneous K\"ahler manifold $X$ is strongly formal over $\Q$.
\end{cor}
\begin{proof}
    Write $X=T\times F$ where $T$ is a torus and $F$ a generalized flag manifold \cite{BorelRemmert}. A product of compact manifolds is strongly formal (over $\Q$) if both factors are strongly formal (over $\Q$)\footnote{See \cite[Prop. 4.21]{MilSt_bigrform}. The argument also works over $\Q$.} and a torus is strongly formal over $\Q$ since it its cohomology algebra is free. On the other hand, any compact K\"ahler manifold satisfies the $\partial\bar\partial$-Lemma and it is known that $H_{dR}(F)$ is an algebra with rational structure of complete intersection type \cite{Borel53HomogeneousCohomologyRegularSequence,BorelHirzebruch58HomogeneousDolbeaultCohomology}.
\end{proof}

\begin{rem}
  A warning on terminology: In \cite{PSZ24} it was shown that compact complex K\"ahler manifolds of dimension $\geq 2$ with the Hodge diamond of a complete intersection are strongly formal. Despite the similarity in terminology, there is very little overlap in these results, as a cohomology algebra with Hodge diamond of a complete intersection is rarely of complete intersection type.
\end{rem}

\section{Proof of \Cref{thm: sfQ}}\label{sec: strong formality over Q}
    Strong formality implies ordinary formality over $\C$, and ordinary formality is independent of the base field. Thus, one might think that any strongly formal cbba is strongly formal over $\Q$. We now show this is not the case. More precisely, we are going to describe a family of cbba's with rational structure $A_\lambda=(A_\Q,A_\C,\varepsilon_\lambda)$ with $\lambda\in \R$ such that $A_\C$ is strongly formal (and in particular also $A_\Q$ is rationally formal), but if $\lambda\in\R\backslash\Q$, $A_\lambda$ is not strongly formal over $\Q$. Just like for smooth complete toric varieties, the cohomology of $A_\lambda$ will be concentrated in diagonal bidegrees.

    \subsection{Definition of the algebras} We set $(A_\Q,d):=(H,0)$, where \[
    H= H(\CP^3\#\CP^3\#(S^2\times S^4);\Q)\]
    and $(A_\C,\partial,\bar\partial):=(H\otimes\C,0,0)$, which we consider as an $\R$-cbba with the bigrading such that $H^{2k}\otimes\C$ is of type $(k,k)$ (note $H^{odd}=0$) and involution $\sigma$ induced the identity on $H$ and conjugation on $\C$. We have $\tot(A_\C)^{\sigma=\id}=H\otimes\R$. Also note that $A_\C$ is tautologically strongly formal.

    \subsection{Definition of $\varepsilon_\lambda$} To define $\varepsilon_\lambda$, first consider the minimal model $\varphi\colon(\Lambda V, d)\rightarrow (H,0)$. Let $\overline{\alpha}\in H^2$ and $\overline{\beta}\in H^4$ denote the generators coming from $S^2$, $S^4$ and let $\overline{x},\overline{y}\in H^2$, be generators coming from the two $\CP^3$ summands such that we have $\overline{\alpha}\overline{\beta} = \overline{x}^3 =\overline{y}^3$.  Generators of $V$ in degrees up to $4$ as well as images in $H$ are given by
    \begin{center}
    \begin{tabular}{c|c|c|c}
        degree & generator & $\xmapsto{d}$ & $\xmapsto{\varphi}$\\
        \hline
          $2$ & $\alpha,x,y$ & $0$ & $\overline{\alpha},\overline{x},\overline{y}$ \\
         $3$ & $a,b,c,e$ & $\alpha^2, x\alpha, y\alpha, xy$ & $0$\\
         $4$ & $\beta,p_1,p_2,p_3,p_4$ & $0,m_1,m_2,m_3,m_4$ & $\overline{\beta},0,0,0,0$
         \end{tabular}
    ´\end{center}
    where
    \begin{align*}
        m_1=xa-\alpha b,\quad m_2=ya-\alpha c,\quad m_3 = yb-\alpha e,\quad m_4 = xc-\alpha e.
    \end{align*}
    $\varphi\colon \Lambda V^{\leq 4}\rightarrow H$ is an isomorphism in cohomology up to degree $5$ and induces a surjection in degree $6$. A basis of the cohomological kernel in degree $6$ is represented by the cocycles
    \begin{center}
    \begin{tabular}{c c c}
        $v_1 = \alpha\beta - x^3$, & $v_2 = \alpha\beta - y^3$, & $v_3 = x\beta$,\\
        $v_4 = y\beta$,& $v_5 = ab + \alpha p_1$,&  $v_6=ac + \alpha p_2$,\\ $v_7 = bc + \alpha p_3 - \alpha p_4 $, & $v_8= ae + xp_2 + \alpha p_4 $,&  $v_9= be + xp_3$, \\
        $v_{10}= ce + y p_4$, & $v_{11} = yp_1 - xp_2 + \alpha p_3 - \alpha p_4~.$ &
    \end{tabular}
    \end{center}
    The computation can be carried out e.g.\ via sage \cite{sagemath}.
    Hence $V^5$ consists of generators $q_i$, $i=1,\ldots,11$, with $d(q_i)=v_i$. We shall not need to compute $V^{\geq 6}$. We set $\xi\colon V^4\rightarrow H^4$ to be the linear map with
    \[\xi(\beta) = 0,\quad \xi(p_1)=-\overline{y}^2,\quad \xi(p_2)=-\overline{x}^2,\quad \xi(p_3)=\overline{\beta}=\xi(p_4).\]
    For any scalar $\lambda\in\R$, we set $\psi_\lambda\colon \Lambda V\otimes\R\rightarrow H\otimes\R$ to be the graded algebra morphism which on $V$ is given by $\varphi + \lambda\xi$ (where $\xi$ is understood as trivial outside of degree $4$). We claim that $\psi_\lambda$ is a cdga morphism. It suffices to check that $\psi_\lambda(dv)=0$ for all generators $v\in V$. For degree reasons this is trivial on $V^{\geq 6}$. Furthermore $\xi$ has no effect on $V^{\leq 3}$ so the condition holds  on $V^{\leq 3}$ due to minimality and $\varphi$ being a cdga map. On $V^4$ the condition is automatic since $H^5=0$. Hence it only remains to check that $\psi_\lambda(dq_i)=\psi_\lambda(v_i)=0$, for $i=1,\ldots,11$, which is indeed the case. Now, $\varepsilon_\lambda$ is the morphism in $\Ho(\cdga_\R)$ represented by the roof

    \[ 
    \begin{tikzcd}
    &\Lambda V\otimes\R\ar[ld, "\psi_\lambda",swap]\ar[rd,"\varphi\,\otimes\,\id_\R"]&\\
    H\otimes\R&&H\otimes\R
     \end{tikzcd}
    \]
    \subsection{Failure of strong formality over $\Q$}
    We first replace $A_\lambda$ by a minimal cbba with rational structure as in \Cref{ex: Mimo for cbba/Q}:
    \begin{lem}
        $A_\lambda$ is quasi-isomorphic to the cbba with rational structure \[A^{min}_\lambda=(\Lambda V,\Lambda W, \tilde{\psi}_\lambda)\]
        where
        \begin{enumerate}
            \item $\Phi:\Lambda W\to A_\C$ is a bigraded minimal model.
            \item $\tilde{\psi}_\lambda$ is a map that makes the following diagram commute in $\Ho(\cdga_\R):$  
               \[
                \begin{tikzcd}
                 &\tot(\Lambda W)^{\sigma=\id}\ar[d,"\Phi"]\\
                \Lambda V\otimes\R\ar[ru, "\tilde{\psi}_\lambda"]\ar[r,"\psi_\lambda"]& H\otimes\R
                \end{tikzcd}
                \] 
        \end{enumerate}
    \end{lem}

  Let us describe $\Lambda W$ and (a choice of) $\tilde{\psi}_\lambda$ explicitly in the degrees relevant to us:
    \begin{center}
    \begin{tabular}{c|c|c|c|c}
       bidegree & generator & $\xmapsto{\del}$  & $\xmapsto{\delbar}$ & $\xmapsto{\Phi}$\\
        \hline
        $(1,1)$ & $\alpha,x,y$ & $0$ & $0$&  $\overline{\alpha},\overline{x},\overline{y}$  \\
         $(1,1)$ & $A,B,C,E$ & $\del A, \del B, \del C, \del E$ & $\delbar A,\delbar B,\delbar C,\delbar E$& $0$\\
        $(2,1)$ & $\del A, \del B,\del C, \del E$ & $0$ & $i\alpha^2, ix\alpha, iy, ixy$ & $0$\\
         $(1,2)$ & $\delbar A, \delbar B,\delbar C, \delbar E$ & $-i\alpha^2, -ix\alpha, -iy, -ixy$ &  $0$ & $0$\\  
          $(2,2)$ & $R_1,\dots,R_4$ & $-i\del m_i'$ &$i\delbar m_i'$& $0$\\
        $(2,1)$ & $P_1,\dots P_4$ & $\del P_i$ &$ R_i-i m_i'$& $0$\\
        $(1,2)$ & $\bar{P}_1,\dots ,\bar{P}_4$ & $R_i + i m_i'$ &$\delbar \bar{P}_i$& $0$\\
        $(3,1)$ & $\del P_1,\dots ,\del P_4$ & $0$ &$\del m_i'$& $0$\\
       
       $(1,3)$ &  $\delbar \bar{P}_1,\dots,\delbar \bar{P}_4$ & $\delbar m_i'$ & $0$& $0$\\
        $(2,2)$ & $\beta$ & $0$ &$0$& $\overline{\beta}$
         \end{tabular}
    \end{center}
    Here, $m_i'$ denote the following decomposable $(2,2)$-classes:
       \begin{align*}
        m_1'=xA-\alpha B,\quad m_2'=yA-\alpha C,\quad m_3'= yB-\alpha E,\quad m_4' = xC-\alpha E.
    \end{align*}

Furthermore, we can take $\tilde{\psi}_\lambda$ as follows (recall $d^c=i(\delbar-\del)$, so $dd^c=2i\del\delbar$):

\begin{center}
    \begin{tabular}{c|c}
    generator & $\xmapsto{\tilde{\psi}_\lambda}$\\\hline
    $\alpha, x,y$&$\alpha, x,y$\\
    $a,b,c,e$& $\frac{1}{2}d^cA, \frac{1}{2}d^cB, \frac{1}{2}d^c C, \frac{1}{2}d^c E$\\
    $\beta$&$\beta$\\
    $p_1,p_2,p_3,p_4$&$\frac{1}{2}R_1 - \lambda y^2, \frac{1}{2}R_2 -\lambda x^2, \frac{1}2 R_3+\lambda\beta, \frac 1 2 R_4 + \lambda \beta$
    \end{tabular}
\end{center}

and we will not need the higher degrees.

    \begin{thm}
        For $\lambda\neq \Q$, the cbba with rational structure $A_\lambda$ is not strongly formal over $\Q$.
    \end{thm}
    
    \begin{proof}
        We apply \Cref{prop: formality of mini cbba/Q}, to $\varphi\colon \Lambda V\rightarrow H$ and $\Phi\colon \Lambda W\rightarrow H\otimes \C$. Identifying cohomologies via these maps, we find that $H(\tilde\psi_\lambda)$ is the identity.
        
        We need to see that there can be no automorphisms $a:\Lambda V\to\Lambda V$ and $b:\Lambda W\to \Lambda W$ of rational cdgas, resp.\ $\R$-cbbas, which induce the identity in cohomology and satisfy $ \Phi \circ b\circ \tilde{\psi}_\lambda\circ (a\otimes \id_\R)=(\varphi\otimes \id_\R) = \psi_0$. Lifting this equality through $\Phi$ we arrive at the condition that \[b\circ \tilde{\psi}_\lambda\circ (a\otimes \id_\R)\simeq \tilde\psi_0.\]
        
        Now, note that 
        \[\pi_{dR}^4(A_\C)=\langle [\beta], [R_1],\dots, [R_4]\rangle,\]
        and 
        \[
        \pi_{dR}^4(A_\C)=\langle [\beta], [p_1],\dots, [p_4]\rangle
        \]
        with induced map $\pi_{dR}^4(\tilde\psi_\lambda)$ sending 
        \[[\beta]\mapsto [\beta],\quad  [p_1],[p_2]\mapsto [\frac 1 2 R_i],\quad  [p_3],[p_4]\mapsto [\frac 1 2 R_i+\lambda\beta]~.
        \]
        As $b$ induces the identity in cohomology it follows that $\pi^4_{dR}(b)[\beta]=[\beta]$. Furthermore 
        \[
        \ker(\pi^4_{dR}(A_\C)\longrightarrow \pi^4_A(A_\C))=\langle [R_1],\dots, [R_4]\rangle
        \]
        needs to be preserved by $\pi^4_{dR}(b)$ since $b$ is an automorphism of cbbas. Hence the coefficient of $[\beta]$ (when expanding an element in the above basis) remains constant under $\pi^4_{dR}(b)$. On the other hand, from the fact that $a$ induces the identity on cohomology we deduce that $a$ preserves all generators of $\Lambda V$ up to degree $3$. As a consequence the $dp_i$ are preserved and it follows that $a(p_i)-p_i$ is closed. In particular it follows that $\pi^4_{dR}(a)[p_i]=[p_i]+\lambda_i[\beta]$ for some $\lambda_i\in \Q$. Thus $\pi^4_{dR}(b\circ \tilde\psi_\lambda\circ a)[p_3]$ has $[\beta]$-coefficient $\lambda-\lambda_3\neq 0$ while the coefficient vanishes for $\pi^4_{dR}(\tilde\psi_0)[p_3]$.
    \end{proof}

    \begin{rem}
        It would be interesting to construct a family of complex manifolds s.t. the associated cbba's with complex structure have a similar behaviour to $A_\lambda$.
    \end{rem}

\subsection{Relation to Mixed Hodge structures} 
In spirit of an argument in \cite{CaCleMo_pi3}, we discuss how the failure of strong formality over $\Q$ can also be detected by a non-trivial extension of Mixed Hodge structures associated with $A_\lambda$. Having given a detailed argument above, we allow ourselves to be brief in describing this second approach:

Let $A=(A_\Q,A_\C,\varphi)$ be a cbba with rational structure such that $A_\C$ satisfies the $\partial\bar\partial$-Lemma. The (total) dual homotopy groups naturally carry rational Mixed Hodge structures. This follows from a discussion analogous to \cite[end of 2.3.1]{StePHT}. 

 To compute this rational Mixed Hodge structure explicitly, we proceed as follows: Let Let $\Lambda V\to A_\Q$ be a minimal model. Consider the filtration $W$ on $\Lambda V$ defined inductively by giving closed generators of degree $k$ weight $k$ letting $d$ and $\wedge$ be strictly compatible with $W_\bullet$. On a minimal model $\Lambda W\to A$ consider the column filtration. The induced filtrations on the (rational structure of) the homotopy bicomplex induce the Mixed Hodge structure on its dual homotopy groups.

We follow the procedure from the previous paragraph in the particular case of $A_\lambda$ (or any $6$-d PD algebra with $H^1=0$). The construction of the minimal model yields a short exact sequence
 \begin{equation}\label{eqn: extn pure Hodge}
0\to C\to \pi^4_{dR}(A_\Q)\to P\to 0~,
 \end{equation}
 where $C$ corresponds to the closed generators in degree $4$ and $P$ can be identified with a complement on which $d$ is injective. By construction the weight of $C$ is $4$. As for the weight of $P$, we observe that $d$ sends $P$ to representatives of triple Massey products of degree $2$ classes. In particular these consist of sums of products of degree $2$ elements (of weight $2$) with degree $3$ elements. The latter are of weight $4$ since $d$ sends them to $\Lambda^2 V^2$ and $V^2\cong H^2$ is of pure weight $2$. Hence $P$ is of weight $6$.
 By the previous discussion, the short exact sequence \eqref{eqn: extn pure Hodge} describes an extension of a weight $6$ Hodge structure by a weight $4$ Hodge structure. Such extensions have been described using intermediate Jacobians in \cite{Carlson_ExtMHS}. Concretely, in the case of $A_\lambda$, using notations as before, we have
\[
C\cong \langle {\beta}\rangle\quad \text{ and }\quad P\cong \langle p_1,..., p_4\rangle\subseteq V^4
\]
and the Hodge structures are pure of type $(2,2)$, resp. $(3,3)$. With these weights and Hodge structures, the description of the group of extensions from \cite{Carlson_ExtMHS} simplifies as follows [loc. cit. p.5]:
\[
\mathrm{Ext}^1_{\Q-MHS}(P,C)\cong \Hom_\C(P,C)/\Hom_\Q(P, C).
\]
The class arising from $A_\lambda$ is given by $[\lambda\cdot \pr\circ \xi|_{P}]$, where $\pr$ denotes the projection to the cokernel of the multiplication $H^2\otimes H^2\to H^4$, which we identify with $C$. In particular the extension is trivial if and only if $\lambda\in \Q$. (We also see that if we consider extensions of real Mixed Hodge structures instead, we get the trivial extension).

\bibliographystyle{alpha}
\bibliography{bib.bib}

\begin{thebibliography}{DGMS75}

\bibitem[AT13]{AnTo}
D.~Angella and A.~Tomassini.
\newblock On the {$\partial\overline{\partial}$}-lemma and {B}ott-{C}hern
  cohomology.
\newblock {\em Invent. Math.}, 192(1):71--81, 2013.

\bibitem[AT15]{AnToBCform}
D.~Angella and A.~Tomassini.
\newblock On {B}ott-{C}hern cohomology and formality.
\newblock {\em J. Geom. Phys.}, 93:52--61, 2015.

\bibitem[BH58]{BorelHirzebruch58HomogeneousDolbeaultCohomology}
A.~Borel and F.~Hirzebruch.
\newblock Characteristic classes and homogeneous spaces. {I}.
\newblock {\em Amer. J. Math.}, 80:458--538, 1958.

\bibitem[Bla56]{Bla_VarAnCplx}
Andr{\'e} Blanchard.
\newblock Sur les vari{\'e}t{\'e}s analytiques complexes.
\newblock {\em Ann. Sci. {\'E}c. Norm. Sup{\'e}r. (3)}, 73:157--202, 1956.

\bibitem[Bod75]{Body_regular}
Richard Body.
\newblock Regular rational homotopy types.
\newblock {\em Comment. Math. Helv.}, 50:89--92, 1975.

\bibitem[Bor53]{Borel53HomogeneousCohomologyRegularSequence}
Armand Borel.
\newblock Sur la cohomologie des espaces fibr\'es principaux et des espaces
  homog\`enes de groupes de {L}ie compacts.
\newblock {\em Ann. of Math. (2)}, 57:115--207, 1953.

\bibitem[BP15]{BuchPan}
Victor~M. Buchstaber and Taras~E. Panov.
\newblock {\em Toric topology}, volume 204 of {\em Math. Surv. Monogr.}
\newblock Providence, RI: American Mathematical Society (AMS), 2015.

\bibitem[BR62]{BorelRemmert}
Armand Borel and Reinhold Remmert.
\newblock {\"U}ber kompakte homogene {K{\"a}hlersche} {Mannigfaltigkeiten}.
\newblock {\em Math. Ann.}, 145:429--439, 1962.

\bibitem[BV83]{BerlineVergne}
Nicole Berline and Mich{\`e}le Vergne.
\newblock {Zeros d’un champ de vecteurs et classes caracteristiques
  equivariantes}.
\newblock {\em Duke Mathematical Journal}, 50(2):539 -- 549, 1983.

\bibitem[Car80]{Carlson_ExtMHS}
James~A. Carlson.
\newblock Extensions of mixed {Hodge} structures.
\newblock Journ\'ees de g\'eometrie alg\'ebrique, {Angers}/{France} 1979,
  107-127 (1980)., 1980.

\bibitem[CCM81]{CaCleMo_pi3}
J.~Carlson, H.~Clemens, and J.~Morgan.
\newblock On the mixed {H}odge structure associated to {$\pi \sb{3}$} of a
  simply connected complex projective manifold.
\newblock {\em Ann. Sci. \'{E}cole Norm. Sup. (4)}, 14(3):323--338, 1981.

\bibitem[CE48]{CheEil_Lie}
Claude Chevalley and Samuel Eilenberg.
\newblock Cohomology theory of {Lie} groups and {Lie} algebras.
\newblock {\em Trans. Am. Math. Soc.}, 63:85--124, 1948.

\bibitem[DGMS75]{DGMS}
P.~Deligne, P.~Griffiths, J.~Morgan, and D.~Sullivan.
\newblock Real homotopy theory of {K}\"ahler manifolds.
\newblock {\em Invent. Math.}, 29(3):245--274, 1975.

\bibitem[DJ91]{DJ}
Michael~W. Davis and Tadeusz Januszkiewicz.
\newblock Convex polytopes, {Coxeter} orbifolds and torus actions.
\newblock {\em Duke Math. J.}, 62(2):417--451, 1991.

\bibitem[Eis95]{Eisenbud_CA}
David Eisenbud.
\newblock {\em Commutative algebra. {With} a view toward algebraic geometry},
  volume 150 of {\em Grad. Texts Math.}
\newblock Berlin: Springer-Verlag, 1995.

\bibitem[FHT01]{FHT_RHT}
Yves Felix, Stephen Halperin, and Jean-Claude Thomas.
\newblock {\em Rational homotopy theory}, volume 205 of {\em Grad. Texts Math.}
\newblock New York, NY: Springer, 2001.

\bibitem[GM13]{GM_RHT}
Phillip~A. Griffiths and John Morgan.
\newblock {\em Rational homotopy theory and differential forms}, volume~16 of
  {\em Prog. Math.}
\newblock New York, NY: Birkh{\"a}user/Springer, 2nd revised and corrected ed.
  edition, 2013.

\bibitem[GS99]{GuilleminSternberg}
Victor~W. Guillemin and Shlomo Sternberg.
\newblock {\em Supersymmetry and equivariant de {R}ham theory}.
\newblock Mathematics Past and Present. Springer-Verlag, Berlin, 1999.
\newblock With an appendix containing two reprints by Henri Cartan [MR0042426
  (13,107e); MR0042427 (13,107f)].

\bibitem[Lil03]{Lillywhite}
Steven Lillywhite.
\newblock Formality in an equivariant setting.
\newblock {\em Trans. Amer. Math. Soc.}, 355(7):2771--2793, 2003.

\bibitem[MP03]{MasudaPanov}
Mikiya Masuda and Taras Panov.
\newblock On the cohomology of torus manifolds.
\newblock {\em Osaka Journal of Mathematics}, 43:711--746, 2003.

\bibitem[MS24]{MilSt_bigrform}
Aleksandar Milivojevi\'c and Jonas Stelzig.
\newblock Bigraded notions of formality and {Aeppli}-{Bott}-{Chern}-{Massey}
  products.
\newblock {\em Commun. Anal. Geom.}, 32(10):2901--2933, 2024.

\bibitem[PR08]{PanovRay}
Taras~E. Panov and Nigel Ray.
\newblock Categorical aspects of toric topology.
\newblock In {\em Toric topology. International conference, Osaka, Japan, May
  28--June 3, 2006}, pages 293--322. Providence, RI: American Mathematical
  Society (AMS), 2008.

\bibitem[PSZ24]{PSZ24}
G.~Placini, J.~Stelzig, and L.~Zoller.
\newblock Nontrivial {M}assey products on compact {K}\"ahler manifolds.
\newblock {\em Preprint arXiv:2404.09867}, 2024.

\bibitem[ST22]{SfTo_DBC}
T.~Sferruzza and A.~Tomassini.
\newblock Dolbeault and {Bott-Chern} formalities: deformations and
  $\partial\bar\partial$-lemma.
\newblock {\em Journal Geom. Phys.}, 175:104470, 2022.

\bibitem[ST24]{SfToBCK}
Tommaso Sferruzza and Adriano Tomassini.
\newblock Bott-{Chern} formality and {Massey} products on strong {K{\"a}hler}
  with torsion and {K{\"a}hler} solvmanifolds.
\newblock {\em J. Geom. Anal.}, 34(11):42, 2024.
\newblock Id/No 348.

\bibitem[Ste21]{StStrDbl}
Jonas Stelzig.
\newblock On the structure of double complexes.
\newblock {\em J. Lond. Math. Soc.}, 104(2):956--988, 2021.

\bibitem[Ste25]{StePHT}
Jonas Stelzig.
\newblock Pluripotential homotopy theory.
\newblock {\em Advances in Mathematics}, 460:110038, 2025.

\bibitem[Sul77]{Sullivan}
D.~Sullivan.
\newblock Infinitesimal computations in topology.
\newblock {\em Inst. Hautes \'{E}tudes Sci. Publ. Math.}, (47):269--331 (1978),
  1977.

\bibitem[{The}25]{sagemath}
{The Sage Developers}.
\newblock {\em {S}ageMath, the {S}age {M}athematics {S}oftware {S}ystem
  ({V}ersion 10.6)}, 2025.
\newblock {\tt https://www.sagemath.org}.

\end{thebibliography}
\end{document}